\documentstyle[12pt]{article}

\begin{document}
\newcommand{\nt}{\noindent}

\begin{center}

{\bf Some Applications of Localization to Enumerative Problems}

\bigskip

Aaron Bertram\footnote{Supported in part by NSF Research Grant
DMS-9970412}

\medskip

Dedicated to Bill Fulton on the occasion of his 60th birthday
\end{center}

\nt {\bf 1. Introduction.} A problem in enumerative geometry 
frequently boils down to the computation of an integral on
a moduli space. We have intersection theory (with Fulton's
wonderful {\it Intersection Theory} \cite{F} as a prime reference) to thank for allowing us to make rigorous
sense of such integrals, but for their computations we
often need to look elsewhere.
If a torus lurks in the background, acting on the moduli space, then
the Atiyah-Bott localization theorem  
allows one to express equivariant cohomology classes on the moduli space 
in terms of their
``residues'' living  on the connected components of the locus of fixed
points (i.e. the fixed submanifolds). This can 
be very useful for computations, particularly when the fixed 
submanifolds are points.

\medskip

We will use localization in a different way.
Here, the moduli space itself will be a fixed 
submanifold for a torus action on a larger ambient space. Localization
is applied in this context to 
relate {\it residues} on the moduli space to 
residues on simpler spaces by means of suitable equivariant 
maps.  This point of view can lead immediately to  
remarkably simple derivations of some complicated-looking formulas,
for example when applied to:

\medskip

(a) Schubert calculus on the flag manifold, and

\medskip

(b) Gromov-Witten invariants of rational curves.

\medskip

In (a) the partial flag manifold
Fl$(1,2,...,m,n)$ is realized as a fixed submanifold of a blown-up
projective space ${\bf P}(\mbox{Hom}(W,V))$ where $W$ and $V$ are
vector spaces of ranks $m$ and $n$ respectively, and all torus actions
come from the ``standard'' torus action of $({\bf C}^*)^m$ on $W^*$. 
The full locus of fixed
points is a disjoint union of
$m!$ fixed submanifolds in this setting, each isomorphic to 
the partial flag manifold, but with different (equivariant) Euler classes. 
For this warm-up application, we will simply list the results of
Jian Kong \cite{K}, where residues on the flag manifold are computed, 
resulting in particular in some new methods for computing Schubert calculus on
the Grassmannian $G(m,n)$. It would be quite interesting to compare this
with other methods (e.g. Gr\"obner bases) for making such computations.

\medskip

In (b), the Kontsevich-Manin moduli space
of stable maps $\overline
M_{0,m}(X,\beta)$ of rational curves 
with image homology class $\beta$ is 
realized as a
fixed submanifold of the ``graph space''
$\overline G_{0,m}(X,\beta) := 
\overline M_{0,0}(X\times ({\bf P}^1)^m,(\beta,1^m))$ again with a
standard torus action. The main applications take place in this setting.

\medskip

When $m = 1$ we investigate
the $J$-functions introduced by Givental in his generalization
of the enumerative side of mirror symmetry to arbitrary
projective manifolds (see \cite{Kim}). The $J$-function is a polynomial 
associated to a complex projective variety $X$ and ample
system of nef divisors which encodes all the one-point Gromov-Witten
invariants. The coefficients of the
$J$-function are push-forwards of residues, and 
our point of view on residues leads to a simple proof
of the multiplicativity of the $J$-functions.  Our point of view
also leads to a non-obvious property of the $J$-function under push-forward.
The $J$-function of projective space is computed in this context
as an immediate consequence of the existence of a nice ``linear'' 
approximation
to the graph space. Following Givental's proof 
of the enumerative mirror conjecture for complete intersections in 
toric varieties, Kim was led to the formulation of a 
``quantum Lefschetz principle'' relating the $J$-function for $X$ with 
$J$-functions for very ample divisor classes in $X$ \cite{Kim}. This has
recently been proved by Y.P. Lee \cite{Lee} in the general case 
building on the proof
in \cite{B} of the case $X = {\bf P}^n$, which we briefly discuss here.

\medskip

When $m > 1$ there are many other fixed submanifolds in the graph 
space besides $\overline M_{0,m}(X,\beta)$, but they all are built
out of Kontsevich-Manin spaces
involving smaller $m$'s and/or smaller $\beta$'s. This 
has been exploited in joint work with Holger Kley 
\cite{BK} to produce recursive formulas for
$m$-point Gromov-Witten invariants, and in particular to 
prove that when the cohomology is generated by divisor classes, the 
$m$-point Gromov-Witten invariants 
can be ``reconstructed'' from one-point Gromov-Witten invariants. 
We will give the formula and an outline of the proof of reconstruction
in the two-point case
as a final application of localization. Another proof
of reconstruction has been
achieved with very different techniques and different formulas
by Lee and Pandharipande \cite{LP}.
As a direct consequence of reconstruction, the
small quantum cohomology of Fano complete intersections in ${\bf P}^n$,
or indeed any toric variety, 
can be explicitly computed, 
since the one-point invariants are computed from
the quantum Lefschetz principle. As another consequence, the 
quantum cohomology of products 
are determined by
reconstruction, since the $J$-functions multiply.

\newpage

\nt {\bf 2. Localization.} When a torus $T = ({\bf C}^*)^m$ acts on a compact
complex  manifold $M$, the fixed submanifolds $F \subset M$
are closed and embedded (of varying dimensions). There is
an equivariant cohomology space H$^*_T(M,{\bf Q})$ which is naturally a module
over the cohomology of the classifying space 
H$^*(BT,{\bf Q}) \cong {\bf Q}[t_1,...,t_m]$. 
If $E$ is a linearized vector bundle over $M$,
then there are equivariant Chern classes $c^T_d(E)$ taking values
in H$^{2d}_T(M,{\bf Q})$, and in particular, the normal bundles
$N_{F/M}$ to the fixed loci are canonically linearized (for the trivial
action of $T$ on $F$) and yield equivariant Euler classes:
$$\epsilon_T(F/M) \in \mbox{H}^*(F,{\bf Q})\otimes_{{\bf Q}}{\bf Q}[t_1,...,t_m]
\cong \mbox{H}^*_T(F,{\bf Q})$$
which are the top equivariant Chern classes of the normal bundles.

\medskip

The
Atiyah-Bott localization theorem \cite{AB} states that these Euler classes are
invertible in H$^*(F,{\bf Q})\otimes_{{\bf Q}}{\bf Q}(t_1,...,t_m)$ and one can
recover an equivariant Chern class $\alpha \in \mbox{H}^*_T(M,{\bf Q})$ uniquely
(modulo torsion) as a sum of residues:
$$\sum_F i_*\frac{i^*\alpha}{\epsilon_T(F/M)}$$
where $i^*$ and $i_*$ are the equivariant pull-back and push-forward
associated to the equivariant inclusion $i:F \hookrightarrow M$. It follows
from the uniqueness that taking residues is functorial. That is, if
$\Phi:M \rightarrow M'$ is an equivariant map and $j:F'\hookrightarrow M'$ is
the inclusion of a component of the fixed submanifold, then:
$$(\dag)\ \ \sum_{F \subset \Phi^{-1}(F')} {\Phi|_F}_*\frac{i^*\alpha}{\epsilon_T(F/M)} = 
\frac{j^*f_*\alpha}{\epsilon_T(F'/M')}$$
where the sum is over the components $X$ of the fixed locus that
are contained in $\Phi^{-1}(F')$ and $\alpha$ is any equivariant
cohomology class on $M$ (see \cite{B} or \cite{LLY}).

\medskip

Thus if we are asked to integrate a cohomology class $\gamma$ on a compact
complex manifold $F$, and if $F$ happens to be isomorphic to a 
component of the fixed locus of an action of $T$ on $M$ as above, then  
the formula above expresses residues at $F$ in terms of residues at $F'$
and at the other fixed loci contained in $\Phi^{-1}(F')$. If $\gamma$ can be 
expressed in terms of residues of equivariant cohomology classes, 
then this formula yields a relation among integrals of cohomology classes
related to $\gamma$. This will be our point of view throughout the rest of
this paper.

\newpage

\nt {\bf 3. Flag Manifolds and Grassmannians.} The partial flag manifold:
$$\mbox{Fl}(1,2,...,m,V) = \{ V_1 \subset V_2 \subset ... \subset V_{m} \subset V \cong {\bf C}^n\ 
|\ V_r \cong {\bf C}^r\}$$  
is a component of the fixed-point locus of an action of $T$ on $M$.

\medskip

In this case, $M$ is the blow-up of ${\bf P}(\mbox{Hom}(W,V))$ 
along:
$$Z_1 \cong {\bf P}(W^*)\times {\bf P}(V) \subset Z_2 \subset ... \subset
Z_{m-1} \subset {\bf P}(\mbox{Hom}(W,V))$$
where $W \cong {\bf C}^m$ and 
$Z_r$ is the locus of maps of rank $\le r$. That is, $M$ is obtained
by blowing up along $Z_1$, followed by the proper transform of $Z_2$,
followed by the proper transform of $Z_3$, etc. If we choose a basis 
$e_1,...,e_m$ of $W$ and let 
$T$ act on the dual space $W^*$with weights $(t_1,...,t_m)$,
then this induces an action of $T$ on $M$, and 
the following are checked in \cite{K}:

\medskip

$\bullet$ The intersection of the $m-1$ exceptional divisors on $M$ is:
$$E_1\cap ... \cap E_{m-1} \cong 
 \mbox{Fl}(1,2,...,W^*) \times \mbox{Fl}(1,2,...,m,V) $$

$\bullet$ The fixed-point loci for the action of $T$ on $M$ are all
contained in this intersection and correspond via the isomorphism above to:
$$ \Lambda_I \times \mbox{Fl}(1,2,...,m,V) $$
where $\Lambda_I$ are the (isolated) fixed points of the action of $T$ on 
Fl$(1,2,...,W^*)$, indexed by the permutations of $m$ letters so that
the permutation $(i_1,...,i_m)$ corresponds to the flag:
$$\Lambda_{(i_1,...,i_m)} = 
\{\left<x_{i_1}\right> \subset \left<x_{i_1},x_{i
_2}\right> \subset \cdots\}$$

\medskip

$\bullet$ Let $\zeta_i$ be the relative hyperplane class for the projection:
$$\mbox{Fl}(1,2,...,i+1,V) \rightarrow \mbox{Fl}(1,2,...,i,V)$$
pulled back to Fl$(1,...,m,V)$. 
Then the equivariant Euler class to the fixed locus 
$F_I = \Lambda_I \times \mbox{Fl}(1,...,m,V) $
is:
$$\epsilon_T(F_I/M) = \prod_{1 \le j < k \le m}(t_{i_k} - t_{i_j})
\prod_{s=1}^{m-1}(t_{i_{s+1}} - t_{i_s} - \zeta_s)$$

We are therefore in a position to apply the formula $(\dag)$ to the diagram:

$$\begin{array}{ccc} M & \stackrel \Phi\rightarrow & M' = {\bf P}(\mbox{Hom}(W,V)) \\
\uparrow \\ F_I \\ \end{array}$$

On the right side, each fixed locus belongs to 
$Z_1 \subset {\bf P}(\mbox{Hom}(W,V))$ as  $F'_i = p_i \times {\bf P}(V)$
where $p_i \in {\bf P}(W^*)$ is the fixed point $\left<x_i\right>$.
In that case, one computes:
$$\epsilon_T(F'_i/M') = \prod_{s \ne i}(h + t_s - t_i)^{n}$$
where $h$ is the hyperplane class on ${\bf P}(V)$. 

\medskip

An $F_I$ belongs to $\Phi^{-1}(F'_i)$ exactly when
$I$ is of the form
$(i,i_2,....,i_m)$. In that case, the induced map:
$$\Phi|_{F_I}: \Lambda_I \times \mbox{Fl}(1,...,m,V) \rightarrow p_i \times {\bf P}(V)$$ 
is the natural projection, which we will denote by $\pi$. Thus $(\dag)$ with $\alpha = 1$ gives us the following interesting 
formula for Schubert calculus:

\medskip

\nt {\bf Schubert Formula 1:}
$$\sum_{\{I|i_1 = i\}} \pi_*\left( \frac 1{\prod_{1 \le j < k \le m}(t_{i_k} - t_{i_j})
\prod_{s=1}^{m-1}(t_{i_{s+1}} - t_{i_s} - \zeta_s)}\right)
= \frac 1{\prod_{s \ne i}(h + t_s - t_i)^{n}}$$

This formula encodes all the information about intersection numbers on the flag
manifold of the form:
$$\int_{{\rm Fl}(1,2,...,m,V)} h^a\cup \zeta_1^{a_1} \cup ...\cup
\zeta_{m-1}^{a_{m-1}}$$

Of course, the same intersection numbers could be obtained by 
applying the Grothendieck relation to each of powers of the $\zeta_i$.
But there is a second formula which is much more interesting, involving
cohomology classes pulled back from the Grassmannian under:
$$\rho:\mbox{Fl}(1,2,...,m,V) \rightarrow G(m,V)$$

\medskip

Recall that such a cohomology class is a symmetric polynomial:
$$\tau(q_1,...,q_m)$$
in the Chern roots $-q_i$ of the universal subbundle $S \subset
V\otimes {\cal O}_{G(m,n)}$.

\newpage

The main theorem of Kong's thesis \cite{K} is the following:

\medskip

\nt {\bf Schubert Formula 2:} 
$$\sum_{\{I|i_1 = i\}} \pi_*\left( \frac {\rho^*\tau(q_1,...,q_m)}{\prod_{1 \le j < k \le m}(t_{i_k} - t_{i_j})
\prod_{s=1}^{m-1}(t_{i_{s+1}} - t_{i_s} - \zeta_s)}\right)$$
$$= \frac {\tau(h + t_1-t_i,....,h + t_m-t_i)}
{\prod_{s \ne i}(h + t_s - t_i)^{n}} + \ \mbox{irrelevant terms}$$
where the irrelevant terms are monomials in the $t_i$ which do 
not appear on the left side of the equation. 

\medskip

\nt {\bf Example:} When $m = 2$, set $i = 1$ above, $t_1 = 0$ and $t_2 = t$. Then:
$$\pi_*\frac {\rho^*\tau(q_1,q_2)}{t(t-\zeta_1)} = \frac {\tau(h,h+t)}{(h+t)^{n}}
+ \mbox{irrelevant terms}$$

If we consider the coefficients of $t^{-2}$ on both sides and integrate,
we get the following new way of doing Schubert calculus on G$(2,V)$:

$$\int_{G(2,V)} \tau(q_1,q_2) = 
\int_{{\rm Fl}(1,2,V)} \pi^*h\cup \rho^*\tau = 
\mbox{coeff of $h^nt^{-2}$ in} \frac{h\cdot\tau(h,h+t)}{(h+t)^{n}}$$ 

\bigskip

Kong proves this formula by finding a suitable equivariant class 
$\alpha$ on $M$ which restricts to the given $\tau$ on each of the 
fixed components $F_I$. This $\tau$ is well enough approximated by 
the pull-back of the corresponding equivariant class of 
a split bundle on $M'$ to give the formula.

\bigskip

The example above for $m=2$ can be similarly worked out for $m > 2$ with 
the main difference being that there are $(m-1)!$ terms on the left
which sum together to the attractive formula on the right. It can 
be shown that this suffices to compute Schubert calculus, and 
it seems that an analysis of the complexity of this computation ought
to be done.

\medskip

Finally, there is no obstruction to carrying out this program when
$V$ is replaced by a vector bundle over a base variety $X$. Kong also
shows how the Chern classes of $V$ figure into this ``relative'' setting
in \cite{K}.

\newpage 

\nt {\bf 4. Gromov-Witten Invariants of Rational Curves.} 
We will describe the relevant Konstsevich-Manin
spaces (and maps among them) only set-theoretically for simplicity.
The interested reader may go to the literature (e.g. \cite{FP}) for 
rigorous constructions of the spaces and morphisms.

\medskip

A map $f:C \rightarrow X$ from an $m$-pointed rational curve is {\bf stable} if:

\medskip

\nt $\bullet$ $C$ has only nodes as singularities, and the marked points are smooth.

\medskip

\nt $\bullet$ Every component of $C$ collapsed by $f$ has at least $3$ 
distinguished points, i.e. marked points and/or nodes.

\medskip

$$\overline M_{0,m}(X,\beta)$$
is the Kontsevich-Manin moduli space of isomorphism classes of stable maps with $m$ marked
points and image homology class $\beta$. If $X$ is ``convex'' 
(e.g. a homogeneous space) then this moduli space is smooth as an orbifold,
of the expected dimension. Otherwise, there is a ``virtual class'' on $X$
with ``all the expected properties'' (see \cite{BF}). There
is always an injective morphism: 
$$\overline M_{0,m}(X,\beta) \hookrightarrow 
\overline G_{0,m}(X,\beta) = \overline M_{0,0}(X\times {({\bf P}^1)}^m,(\beta,1^m))$$
where the latter space is the ``graph space'' associated to the former.

\medskip

Given a stable map $f:C \rightarrow X$ and points $p_1,...,p_m\in C$,
we obtain the image of $[f]$ in the graph space by attaching a copy
of ${\bf P}^1$ to each of the points, gluing $p_i\in C$ to $0\in {\bf P}^1$,
and collapsing each ${\bf P}^1$ to construct the resulting stable map 
$g: C \cup \coprod {\bf P}^1 \rightarrow X$. 

\medskip

It is convenient to number the ${\bf P}^1$'s, so ${\bf P}^1_i = {\bf P}(W_i)$
is the particular ${\bf P}^1$ which we attach to $p_i$. The
actions of ${\bf C}^*$ on the dual spaces $W_i^*$ with weights $(0,t_i)$
give a natural action of the torus $T$ on the product of the ${\bf P}^1$'s and hence
on the graph space above. Moreover, the $m$-pointed Kontsevich-Manin
space is one of the components of the fixed-locus for the torus 
action. 

\medskip

One computes, using for example \cite{FP}:
$$\epsilon_T(\overline M_{0,m}(X,\beta)/\overline G_{0,m}(X,\beta)) = 
\prod_{i=1}^mt_i(t_i - \psi_i)$$
where the $\psi_i$ are the ``gravitational descendants''
$\psi_i = c_1(\sigma_i^*(\omega))$. Here $\omega$ is the relative dualizing sheaf
of the universal curve ${\cal C}$ over $\overline M_{0,m}(X,\beta)$ and $\sigma_i$
is the section of ${\cal C}$ corresponding to the $i$th marked point.

\bigskip

\nt {\bf The Case $m=1$:}  Here we let $t = t_1$ and $\psi = \psi_1$.

\medskip

If $H$ is an ample divisor on $X$, then following Givental, we define:
$$J_{X,H}(q) = 1 + \sum_{\beta \ne 0} {e_\beta}_*\left(\frac 1{t(t-\psi)}\right) q^{\rm{deg}_H(\beta)}$$
where $e_\beta:\overline M_{0,1}(X,\beta) \rightarrow X$ is the evaluation map 
$e_\beta([f]) = f(p)$. Since only a finite number of classes $\beta$ have a 
given degree against $H$, this sum makes sense. More generally, we will suppose
$H$ is a system $H = (H_1,...,H_r)$ of (linearly independent) nef divisors, and
that some linear combination of the $H_i$ is ample. In that case, we define:
$$J_{X,H}(q) = 1 + \sum_{\beta \ne 0} {e_\beta}_*\left(\frac 1{t(t-\psi)}\right) 
q_1^{{\rm deg}_{H_1}(\beta)}...q_r^{{\rm deg}_{H_r}(\beta)}$$

\medskip

The following ``functorial'' properties of the $J$-function are easily
proved
once we recognize that the coefficients are push-forwards of residues. 

\medskip

\nt {\bf Product Formula:} Suppose $X$ and $X'$ are simply connected projective
manifolds (so the curve classes on $X \times X'$ are all of the form $(\beta,\beta')$)
and $H$ and $H'$ are ample systems of divisors, as above. Then:
$$J_{X\times X',(\pi_1^*H,\pi_2^*H')}(q,q') = 
\pi_1^*J_{X,H}(q) \cdot \pi_2^*J_{X',H'}(q')$$

{\bf Proof:} Kontsevich-Manin spaces are functorial, in the sense that a map 
$f:X \rightarrow Y$ gives rise to maps:
$$f_{0,m}: \overline M_{0,m}(X,\beta) \rightarrow 
\overline M_{0,m}(Y,f_*\beta)$$
and analogous compatible equivariant maps on the graph spaces. 
Thus the projection maps give rise to a 
diagram of ``lifts'' of the identity map:
$$\begin{array}{ccc} \overline G_{0,1}(X\times X',(\beta,\beta')) & 
\stackrel \Phi\rightarrow &
\overline G_{0,1}(X,\beta) \times \overline G_{0,1}(X',\beta')\\
\uparrow & & \uparrow \\
\overline M_{0,1}(X\times X',(\beta,\beta')) & \stackrel \phi \rightarrow &
\overline M_{0,1}(X,\beta) \times \overline M_{0,1}(X',\beta')\\  
\downarrow & & \downarrow \\
X\times X' & \rightarrow & X\times X'
\end{array}$$

\nt $\Phi$ is birational when $X$ and $X'$ are convex (and ``virtally birational'' always)
even though $\phi$ is not birational (the two sides have
different dimensions!). Thus $\Phi_*1 = 1$ and we may apply $(\dag)$ to the class $1$ to obtain:
$$\phi_*\left(\frac 1{t(t-\psi)}\right) = \pi_1^*\left(\frac 1{t(t-\psi)}\right)
\pi_2^*\left(\frac 1{t(t-\psi)}\right)$$
Further pushing forward to $X\times X'$ yields the desired product formula.

\medskip

\nt {\bf Push-Forward Formula:} Suppose $f:X\rightarrow Y$ is given. Then
there are equivariant classes 
${f_\beta}_*1 \in {\rm H}^*(\overline M_{0,1}(Y,f_*\beta),{\bf Q})[t]$
so that:
$$f_*J_{X,H}(q) = f_*1 + \sum_{\beta \ne 0} {(e_{f_*\beta})}_*\left(\frac {{f_\beta}_*1}
{t(t-\psi)}\right) q_1^{{\rm deg}_{H_1}(\beta)}...q_r^{{\rm deg}_{H_r}(\beta)}$$

{\bf Proof:} Here we consider the diagram of lifts of $f$:
$$\begin{array}{ccc} \overline G_{0,1}(X,\beta) & 
\stackrel \Phi\rightarrow &
\overline G_{0,1}(Y,f_*\beta)\\
\uparrow & & j\uparrow \\
\overline M_{0,1}(X,\beta) & \stackrel \phi \rightarrow &
M_{0,1}(Y,f_*\beta)\\  
\downarrow & & \downarrow \\
X & \stackrel f\rightarrow & Y
\end{array}$$
and note that applying $(\dag)$ to the class $1$ again, we get:
$$f_*{e_\beta}_*\left( \frac 1{t(t-\psi)}\right) =
(e_{f_*\beta})_*\phi_*\left(\frac 1{t(t-\psi)}\right) = 
(e_{f_*\beta})_*\left( \frac{j^*\Phi_*1}{t(t-\psi)}\right)$$
hence the push-forward formula with ${f_\beta}_*1 := j^*\Phi_*1$.

\medskip

\nt {\bf Remark:} If $f$ is an embedding, then $\phi^*\psi = \psi$, in which 
case the projection formula tells us that ${f_{\beta}}_*1 = \phi_*1$ is constant
in $t$. It seems 
that in general, however, ${f_{\beta}}_*1$ is not constant in $t$. It
would be very interesting to compute it, for instance, 
in case $f$ is the inverse
of a blow-up along a submanifold.

\medskip

\nt {\bf The $J$-function of Projective Space:} Let $H$ be the hyperplane class
on ${\bf P}^n$. Then:
$$J_{{\bf P}^n,H}(q) = 
\sum_{d = 0}^\infty {e_d}_*\left(\frac 1{\prod_{k=1}^d(H + kt)^{n+1}}\right)q^d$$

{\bf Proof:} $\overline G_{0,1}({\bf P}(V),d)$ 
has a natural birational map to a ``linear'' space
${\bf P}(\mbox{Hom}(\mbox{Sym}^d(W),V)$, where $W = W_1$. A general element of 
the graph space is represented by a degree $d$ morphism $f:{\bf P}^1 \rightarrow {\bf P}^n$
which maps to an $n+1$-tuple
of degree $d$ polynomials $(p_0(x,y):...:p_n(x,y))$ with no common factors.
When the curve underlying the
stable map picks up extra components, then the $n+1$-tuple of 
polynomials picks up common factors. In particular, the image of 
$\overline M_{0,1}({\bf P}^n,d)$ under this weighted blow-down is
a copy of ${\bf P}^n$, embedded via:
$$\{x^d\} \times {\bf P}(V) \hookrightarrow 
{\bf P}(\mbox{Sym}^d(W^*))\times {\bf P}(V) \hookrightarrow  
{\bf P}(\mbox{Hom}(\mbox{Sym}^d(W),V)$$

Thus, we have the diagram:
$$\begin{array}{ccc} \overline G_{0,1}({\bf P}^n,d) &
\stackrel \Phi\rightarrow & {\bf P}(\mbox{Hom}(\mbox{Sym}^d(W),V) \\ \\
i\uparrow &  & j\uparrow \\ \\
\overline M_{0,1}({\bf P}^n,d) & \stackrel {e_d} \rightarrow & {\bf P}^n
\end{array}$$ 
One computes (see \cite{BDPP})
$\epsilon_T({\bf P}^n/{\bf P}(\mbox{Hom}(\mbox{Sym}^d(W),V))) = 
\prod_{k=1}^d(H+kt)^{n+1}$
so that $(\dag)$ now applies with the class $1$, giving us:
$${e_d}_*\left(\frac 1{t(t-\psi)}\right) = 
\frac{1}{\prod_{k=1}^d (H+kt)^{n+1}}$$
proving the formula.

\medskip

\nt {\bf Quantum Lefschetz Hyperplane:} We will limit ourselves to considering
hypersurfaces in ${\bf P}^n$, as in \cite{B}. See \cite{Lee} for the general version. 
Let $f:X \hookrightarrow {\bf P}^n$ be a hypersurface of degree $l$, and let
$H$ denote the hyperplane class, either on ${\bf P}^n$ or on $X$. Let:

$$I_{X/{\bf P}^n,H}(q) = \sum_{d=0}^\infty \frac{\prod_{k=0}^{dl}(lH + kt)}
{\prod_{k=1}^d(H + kt)} q^d$$ 

(a) If $l < n$, then $f_*J_{X,H}(q) = I_{X/{\bf P}^n,H}(q)$.

\medskip

(b) If $l = n$, then $f_*J_{X,H}(q) = e^{\frac{-l}tq}I_{X/{\bf P}^n,H}(q)$.

\medskip

(c) If $l = n+1$, then there are $a(q),b(q) \in q{\bf Q}[[q]]$ so that:
$$f_*J_{X,H}(q) = e^{\frac Hta(q) + b(q)}I_{X/{\bf P}^n,H}(qa(q))$$

To prove this, one uses the diagram for ${\bf P}^n$ and observes that 
$$f_*J_{X,H}(q) = lH + 
\sum_{d>0} {e_d}_*\left(\frac {j^*\Phi_*[X]_T}{t(t-\psi)}\right)q^d$$
where $[X]_T$ is the equivariant Chern class:
$$[X]_T = c_{dl+1}^T(\pi_*e_d^*{\cal O}_{{\bf P}^n}(l)) \in
\mbox{H}^*_T(\overline G_{0,1}({\bf P}^n,d),{\bf Q})$$

Thus the proof of quantum Lefschetz amounts to a 
detailed analysis of the class
$j^*\Phi_*[X]_T \in H^*({\bf P}^n,{\bf Q})[t]$.
This is obtained by decomposing $[X]_T$ along boundary strata of
the graph space by means of intersection theory. 
In particular, the open stratum of the 
graph space contributes
$\prod_{k=0}^{dl}(lH+kt)$
via an approximation of $\pi_*e_d^*{\cal O}_{{\bf P}^n}(l)$ 
by $\mbox{Sym}^{dl}(W^*) \otimes \Phi^*{\cal O}(l)$
in much the same way that the graph space is approximated by 
${\bf P}(\mbox{Hom}(\mbox{Sym}^{d}(W),V))$. In the case $l < n$,
this is the only stratum which contributes to $j^*\Phi_*[X]_T$,
giving us (a). In the other cases, the boundary strata do 
contribute, but
in a self-similar manner. When tallied up, these contributions 
give formulas (b) and (c) in the cases $l = n$ and $l = n+1$
respectively. It is unknown whether a more general ``change of
coordinates'' analogous to (b) and (c) occurs in the general
type cases $l > n+1$. 

\bigskip

\nt {\bf The case $m > 1$. Reconstruction:} In \cite{BK}, reconstruction
theorems make use of the following 
diagrams of K-M spaces and graph spaces:
$$\begin{array}{ccc}
\overline G_{0,m}(X,\beta) & 
\stackrel \Phi \rightarrow
& \overline G_{0,m-1}(X,\beta) \times 
\overline M_{0,0}({\bf P}^1_1 \times {\bf P}^1_m,(1,1)) 
\\ \\ i\uparrow & & j\uparrow \\ \\
\overline M_{0,m}(X,\beta) & \stackrel \pi \rightarrow & 
\overline M_{0,m-1}(X,\beta)
\end{array}$$ 

$\Phi$ is derived, as in the product formula, from
projections. $\overline M_{0,m-1}(X,\beta)$ is included 
in the graph space in the ordinary way, and the inclusion
$j$ is given by the additional inclusion of the point 
corresponding to the inclusion of 
the intersecting lines $\{0\}\times {\bf P}^1_m \cup 
{\bf P}^1_1 \times \{0\}$ in ${\bf P}^1_1 \times {\bf P}^1_m$.

\medskip

The fixed loci contained in $\Phi^{-1}(\overline M_{0,m-1}(X,\beta))$,
in addition to $\overline M_{0,m}(X,\beta)$ are isomorphic to
one of the following:
$$\overline M_{0,k+1}(X,\beta_1) \times_X 
\overline M_{0,m-k}(X,\beta - \beta_1)) \ 
\mbox{or}\  \overline M_{0,m-1}(X,\beta)$$
and the induced maps to $\overline M_{0,m-1}(X,\beta)$ are 
the gluing maps to boundary divisors (see \cite{FP}) and
the identity map, respectively.  

\medskip

The equation $(\dag)$ now tells us that given an equivariant cohomology
class $\alpha$ on $\overline G_{0,m}(X,\beta)$, there is 
a relation among the residues of $\alpha$ along the fixed
loci listed above, as well as the residue of $\Phi_*\alpha$ along
the fixed locus $\overline M_{0,m-1}(X,\beta)$. So the question
now becomes, how to find interesting equivariant classes $\alpha$
on the graph space? The only source we know of to produce
good residue classes comes from
the linear approximation to $\overline G_{0,1}({\bf P}^n,d)$.
Namely, suppose a morphism (not necessarily an embedding)
$f:X \rightarrow {\bf P}^n$ is given. Then we can pull back equivariant
cohomology classes via:
$$\overline G_{0,m}(X,\beta) \rightarrow \overline G_{0,1}(X,\beta)
\rightarrow \overline G_{0,1}({\bf P}^n,d) \rightarrow 
{\bf P}(\mbox{Hom}(\mbox{Sym}^d(W),V))$$

After all the equivariant Euler classes are computed, recursive
formulas are obtained. Thus in this context the necessity of
considering cohomology classes generated by divisor classes
springs from our inability to find useful equivariant classes
not coming from the linear approximation spaces to 
${\overline G}_{0,1}({\bf P}^n,d)$. 
As an example of the reconstruction theorems we obtain,
we include the most useful one, which, in case the cohomology
of $X$ is generated by divisor classes, already suffices to 
express (small) quantum cohomology in terms of the $J$-function.

\medskip

\nt {\bf Reconstruction Theorem for $2$-Point Invariants:} Given
$f:X \rightarrow {\bf P}^n$, let $H$ be the hyperplane class
on ${\bf P}^n$ and on $X$, and define:

\medskip

$F_\beta(t) = {e_\beta}_*\left( \frac 1{t(t-\psi)}\right)$ (these are the coefficients
of $J$) and

\medskip

$G_\beta(\gamma,t) = {e_\beta^1}_*\left(\frac{{e_\beta^2}^*\gamma}{t-\psi_2}\right)
\ \mbox{for evaluation maps}\ 
e_\beta^1,e_\beta^2:M_{0,2}(X,\beta) \rightarrow X\ \mbox{and}$ 
cohomology class $\gamma \in {\rm H}^*(X,{\bf Q})$,
extended by linearity in the 
first factor. 

\medskip

Then the expression:
$$G_\beta(H^a,t) + 
\left(\sum_{\beta_1 + \beta_2 = \beta} 
G_{\beta_1}(F_{\beta_2}(-t)(H-d_{\beta_2}t)^a,t)\right) + 
F_{\beta}(-t)(H-d_\beta t)^a$$
is polynomial in $t$, where $d_\beta$ is the degree of 
$f_*\beta \in {\rm H}_2({\bf P}^n,{\bf Z})$.

\bigskip

Since $G_\beta(H^a,t)$ is polynomial in 
$t^{-1}$ (with no contant term), this formula expresses
$G_\beta(H^a,t)$ in terms of coefficients of $J$ and $G_{\beta '}(H^a,t)$
for smaller $\beta'$. Hence it inductively determines 
$G_\beta(H^a,t)$
in terms of $J$.

\bigskip

\nt University of Utah, Salt Lake City, UT 84112

\medskip

\nt bertram@math.utah.edu

\end{document}